\documentclass[a4paper,11pt]{article}
\usepackage[scale=0.75,twoside,bindingoffset=5mm]{geometry}

\usepackage[english]{babel}
\usepackage[utf8]{inputenc}
\usepackage{amsmath}
\usepackage{graphicx}
\usepackage{subcaption}
\usepackage{amssymb}
\usepackage{amsthm}
\usepackage{tikz-cd}
\usepackage{mathrsfs}
\usepackage[colorinlistoftodos]{todonotes}
\usepackage{enumitem}
\usepackage{yfonts}
\usepackage{ dsfont }
\usepackage{setspace}
\usepackage{mathtools}

\usepackage[title]{appendix}

\title{Phase Transition in the Generalized Stochastic Block Model}

\author{Sun Min Lee\footnote{Department of Mathematical Sciences, KAIST, Daejeon 34141, Korea
	\newline e-mail: \texttt{adelialee@kaist.ac.kr}} \and	Ji Oon Lee\footnote{Department of Mathematical Sciences, KAIST, Daejeon 34141, Korea
	\newline email: \texttt{jioon.lee@kaist.edu}}}

\date{\today}
\newtheorem{thm}{Theorem}[section]

\numberwithin{equation}{section}

\theoremstyle{definition}
\newtheorem{defn}{Definition}[section]

\theoremstyle{remark}

\newcommand{\R}{\mathbb{R}}
\newcommand{\C}{\mathbb{C}}

\newcommand{\N}{\mathbb{N}}

\newcommand{\E}{\mathbb{E}}

\newcommand{\dd}{\mathrm{d}}
\newcommand{\ii}{\mathrm{i}}

\newcommand{\wt}{\widetilde}
\newcommand{\wh}{\widehat}

\DeclareMathOperator{\Tr}{Tr}
\DeclareMathOperator{\re}{Re}
\DeclareMathOperator{\im}{Im}

\newcommand{\beq}{ \begin{equation} }
\newcommand{\eeq}{ \end{equation} }

\begin{document}
\maketitle

\begin{abstract}
We study the problem of detecting the community structure from the generalized stochastic block model (GSBM). Based on the analysis of the Stieljtes transform of the empirical spectral distribution, we prove a BBP-type transition for the largest eigenvalue of the GSBM. For specific models such as a hidden community model and an unbalanced stochastic model, we provide precise formulas for the two largest eigenvalues, establishing the gap in the BBP-type transition.
\end{abstract}

\section{Introduction} \label{sec:intro}

One of the most fundamental and natural problems in data science is to understand an underlying structure from data sets that can be viewed as networks. The problem is known as the clustering or community detection, and it appears in diverse field of studies involving real-world networks. 

The stochastic block model (SBM) is one of the most fundamental mathematical models to understand the community structure in networks. An SBM is a random graph with $N$ nodes, partitioned into $K$ disjoint subsets, called the communities, $C_1, C_2, \dots, C_K$. One can characterize an SBM via its adjacency matrix, which is a symmetric (random) matrix $\wt M$, whose $(i, j)$-entry $\wt M_{ij}$ is a Bernoulli random variable depending only on the communities to which the nodes $i$ and $j$ belong. For the clustering of an SBM, it is often useful to analyze the eigenvalues of the adjacency matrix and their associated eigenvectors, which is known as a spectral method. 

One of the most prominent examples of spectral methods is the principal component analysis (PCA) in which the behavior of the eigenvectors associated with the extremal eigenvalues are considered to obtain the community structure of the SBM. For an SBM with two communities, the expectation of its adjacency matrix $\wt M$ has a block structure, i.e.,
\beq \label{eq:2by2}
	\E[\wt M] =\left(
	\begin{array}{c|c}
		P_{11} & P_{12} \\
		\hline
		P_{21} & P_{22}
	\end{array}
	\right).
\eeq
In the simplest case of a balanced SBM with $P_{11} = P_{12} = p$, $P_{12} = P_{21} = q$, and the two communities are of equal size, it can be easily checked that $\E[\wt M]$ has at most two non-zero eigenvalues, $N(p+q)/2$ and $N(p-q)/2$. Thus, if $N(p-q)/2$ is sufficiently large, the perturbation $\wt M - \E[\wt M]$ is negligible for the two largest eigenvalues $\wt M$ and it is possible to determine the community structure from the eigenvector associated with the second largest eigenvalue of $\wt M$. Equivalently, after subtracting $(p+q)/2$ from each entry, the (shifted) adjacency matrix becomes the sum of a rank-$1$ deterministic matrix and a random matrix with centered entries, and one can use the eigenvector associated with the largest eigenvalue of the shifted adjacency matrix for clustering.

The sum of a deterministic matrix and a random matrix has been extensively studied in random matrix theory. When the deterministic matrix is rank-$1$ and the random matrix is a Wigner matrix, it is called a (rank-$1$) spiked Wigner matrix. The behavior of the largest eigenvalue of a spiked Wigner matrix is known to exhibit a sharp phase transition depending on the ratio between the spectral norms of the deterministic part and the random part. This type of phase transition is called the BBP transition, after the seminal work of Baik, Ben Arous, and P\'ech\'e \cite{BBP2005} for spiked (complex) Wishart matrix. From the BBP transition, we can immediately see that the detection of the signal is possible via PCA when the signal-to-noise ratio (SNR) is above a certain threshold.

While the BBP transition has been proved for spiked Wigner matrices under various assumptions \cite{Peche2006,FeralPeche2007,CapitaineDonatiFeral2009,Raj2011}, it is not directly applicable to the SBM, since the entries in a Wigner matrix are i.i.d. (up to symmetry constraint) whereas those in the adjacency matrix of an SBM are not. The proof of the BBP transition with an SBM is substantially harder. For example, unless the SBM is balanced, the empirical spectral distribution (ESD) of $\wt M$ does not even converge to the semi-circle distribution, which is the limiting ESD of a Wigner matrix; the limiting ESD in this case is not given by a simple formula as the semi-circle distribution but by an implicit formula via its Stieltjes transform.

\subsection*{Main contribution}

In this paper, we consider a model that generalizes the SBM, called the generalized stochastic block model (GSBM), with two communities. In this model, the mean of the matrix has the same block structure as that of the SBM in \eqref{eq:2by2}, but the entries are not necessarily Bernoulli random variables. See Definition \ref{def:GSBM} for the precise definition of the GSBM. 

For the GSBM, We prove the BBP-type transition for its largest eigenvalue (Theorem \ref{thm:main}). The proof is based on the analysis of the Stiejtes transform of the ESD, which involves the resolvent of the random part of the GSBM. Due to the community structure, the random part is not a Wigner matrix, but a generalization of a Wigner matrix, known as a Wigner-type matrix. The local properties of eigenvalues of Wigner-type matrices are now well-established by recent developments of random matrix theory; see, e.g., \cite{ajanki2019quadratic,AEK17,dumitriu2019sparse}.

In our main result, Theorem \ref{thm:main}, we only state the existence of the critical values and the limiting gap between the two largest eigenvalue but refrain from writing the precise formulas for them. We instead apply our results to specific examples naturally arising in applications, hidden community model and unbalanced stochastic block model, and present the results from numerical experiments. (In terms of the edge probability, the former corresponds to the case $P_{11} = p$ and $P_{12} = P_{21} = P_{22} = q$, while the latter $P_{11} = P_{22} = p$ and $P_{12} = P_{21} = q$ (but $\gamma \neq 1/2$)).

\subsection*{Related works}
The local law for Wigner-type matrices and the behavior of Quadratic vector equations (QVE), which are crucial in the analysis for Wigner-type matrices, were thoroughly investigated by Ajanki, Erd{\H{o}}s and Kr{\"u}ger \cite{ajanki2019quadratic,AEK17}. A related result on the local law at the cusp for the Wigner-type matrix was also proved \cite{erdHos2020cusp}. For more results on general Wigner-type matrices, we refer to \cite{dumitriu2019sparse, erdHos2019bounds, zhu2020graphon} and references therein.

The phase transition of the largest eigenvalue was first proved proved by Baik, Ben Arous and P\'ech\'e \cite{BBP2005} for spiked Wishart matrices and later extend to other models, including the spiked Wigner matrix under various assumptions \cite{Peche2006,FeralPeche2007,CapitaineDonatiFeral2009,Raj2011}. If the SNR is below the threshold given by the BBP transition, the largest eigenvalue has no information on the signal and we cannot use the PCA for the detection of the signal. For this case, the PCA can be improved by an entrywise transformation that effectively increase the SNR \cite{NIPS_Barbier,Perry2018}. Reliable detection is impossible below a certain threshold \cite{Perry2018}, and it is only possible to consider a weak detection, which is a hypothesis testing between the null model (without spike) and the alternative (with spike). For more detail about the weak detection, we refer to \cite{AlaouiJordan2018,chung2019,jung2020}.

The problem of recovering a hidden community from a symmetric matrix for two important cases, the Bernoulli and Gaussian entries, was discussed by Hajek, Wu, and Xu \cite{Hajek-2017}. A threshold for exact recovery in SBM was discussed in \cite{AbbeEmm2017,AbbeEmm2014,Chen-2015,hajek2016achieving}. Recovering community at the Kesten--Stigum threshold for SBM was considered in \cite{hajek2018recovering}. For more results and Applications on SBM, we refer to \cite{stanley2019stochastic} and references therein.

\subsection*{Organization of the paper}

The rest of the paper is organized as follows: In Section \ref{sec:main}, we define the model and state the main result. In Section \ref{sec:ex}, we introduce the hidden community model and unbalanced stochastic model to provide the results from numerical experiments around the transition threshold. In Section \ref{sec:sketch}, we prove the main theorem. A summary of our results and future research directions was discussed in Section \ref{sec:conclusion}. Appendix \ref{sec:prelim} contains the definition of the Wigner-type matrices and preliminary results on this model. The detailed analysis for the specific models can be found in Appendix \ref{sec:sbm}.

\section{Main Results} \label{sec:main}

In this section, we precisely define the matrix model that we consider in this paper and state our main theorem. We begin by introducing a shifted, rescaled matrix for a generalized stochastic block model with two communities.

\begin{defn}[Generalized Stochastic Block Model (GSBM)] \label{def:GSBM}
An $N \times N$ matrix $M$ is a generalized stochastic block model if
\[
	M = H + \lambda uu^T
\]
where $\lambda \geq 0$ is a constant, $u = (u_1, u_2, \dots, u_N) \in \R^N$ with $\|u\| = 1$, and $H=[H_{ij}]$ is an $N\times N$ real symmetric matrix, satisfying the following:

\begin{itemize}

\item There exist $S \subset [N] := \{1, 2, \dots, N \}$ and constants $\theta_1, \theta_2$ such that
		\[
			u_i = 
				\begin{cases}
					\theta_1 & \text{if} \quad i \in S \,, \\
					\theta_2 & \text{if} \quad i \notin S \,.
				\end{cases}
		\]
	We further assume that $\frac{|S|}{N}, (1-\frac{|S|}{N}) > c >0$ for some ($N$-independent) constant $c$.

\item Upper diagonal entries $H_{ij} (i\le j)$ are centered independent random variables such that

		\begin{itemize}
		\item there exist ($N$-independent) constants $\alpha_1$ and $\alpha_2$ such that
		\[
			\E[H_{ij}^2] = 
				\begin{cases}
					\alpha_1 N^{-1} & \text{if} \quad i, j \in S \\
					\alpha_2 N^{-1} & \text{if} \quad i, j \notin S \\
					N^{-1} & \text{otherwise}
				\end{cases}
		\]
		\item for any ($N$-independent) $D>0$, there exists a constant $C_D$ such that for all $i \leq j$
		\[
			\E[H_{ij}^D] \leq C_D N^{-\frac{D}{2}}.
		\]
		\end{itemize}
\end{itemize}
\end{defn}

For an adjacency matrix $\wt M$ in \eqref{eq:2by2}, if $P_{11} = p_1$, $P_{22} = p_2$, and $P_{12} = P_{21} = q$, then after shifting and rescaling, we find that
\beq \label{eq:p_12}
	\alpha_1=\frac{p_1(1-p_1)}{q(1-q)}, \qquad \alpha_2=\frac{p_2(1-p_2)}{q(1-q)}.
\eeq
(See Appendix \ref{sec:sbm} for more detail.)

We remark that $H_{ij}$ are not necessarily Bernoulli random variables. The assumption on the finite moment means that the model is in the dense regime.
The most typical balanced stochastic block model with two communities correspond to the choice of parameters $|S| = N/2$ and $\alpha_1 = \alpha_2 > 1$.

Our main theorem is the following result on the phase transition for the spectral gap of GSBM.

\begin{thm} \label{thm:main}
Let $M$ be a generalized stochastic block model defined in Definition \ref{def:GSBM}. Denote by $\lambda_1$ and $\lambda_2$ the largest and the second largest eigenvalue of $M$. Assume that $\gamma := N_1 / N$ is fixed. Then, there exists a constant $\lambda_c$, depending only on $\theta_1, \theta_2, \alpha_1, \alpha_2, \gamma$, such that
\begin{itemize}
\item (Subcritical case) if $\lambda < \lambda_c$, then $\lambda_1 - \lambda_2 \to 0$ as $N \to \infty$, almost surely.
\item (Supercritical case) if $\lambda > \lambda_c$, then $\lambda_1 - \lambda_2 \to g$ as $N \to \infty$, almost surely, for some ($N$-independent) constant $g$.
\end{itemize}
\end{thm}

We do not include the precise formulas for the critical value $\lambda_c$ and the gap $g$ in the statement of Theorem \ref{thm:main} for the general cases since they are lengthy but not particularly informative. The formulas for the special cases can be found in Section \ref{sec:ex}.

\section{Examples and Experiments} \label{sec:ex}

In this section, we focus on several specific models and check how the main result, Theorem \ref{thm:main}, applies to them. 

\subsection{Hidden community model}\label{sec:exh}

In the hidden community model, only one of the intra-community connection probability is larger than the inter-community connection probability, and the other intra-connection probability coincides with the inter-community connection. The precise definition for such a model is as follows:

\begin{defn}[Hidden Community Model] \label{defn:hidden} Let $C\subset[n]$ such that $|C|=K$. Let define that $S$ is a $N\times N$ symmetric matrix with $S_{ii}=0$ where $S_{ij}$ are independent for $1 \leq i \leq j \leq N$ and
\begin{align*}
S_{ij}\sim
 \begin{cases}
    P, & \mbox{if }i,j\in C \\
    Q, & \mbox{otherwise} \\
 \end{cases}
\end{align*}
for given probability measures $P$ and $Q$. 
\end{defn}

We consider the BBP-type transition of the hidden community model with Bernoulli entries, i.e., $P=\mathrm{Bernoulli}(p)$ and $Q=\mathrm{Bernoulli}(q)$ with $p\neq q$, which also corresponds to the case $\alpha_2 = 1$ or $p_2 = q$ in \eqref{eq:p_12}. It is not hard to find that the transition occurs in the regime
\[
	p_1 := p = \frac{w}{\sqrt{N}}+q
\]
for some (possibly $N$-dependent) $w = \Theta(1)$. After shifting and rescaling, we find that $\lambda_2 \to 2$ and
\beq
	\lambda_1 \to
		\begin{cases}
			\frac{\gamma w}{\sqrt{q(1-q)}} + \frac{\sqrt{q(1-q)}}{\gamma w} & \text{ if } w > \frac{\sqrt{q(1-q)}}{\gamma} \,, \\
			2 & \text{ if } w < \frac{\sqrt{q(1-q)}}{\gamma} \,.
		\end{cases}
\eeq
See Appendix \ref{subsec:hidden} for the detail.

\begin{figure}[h]
\centering
\subfloat{{\includegraphics[width=7cm]{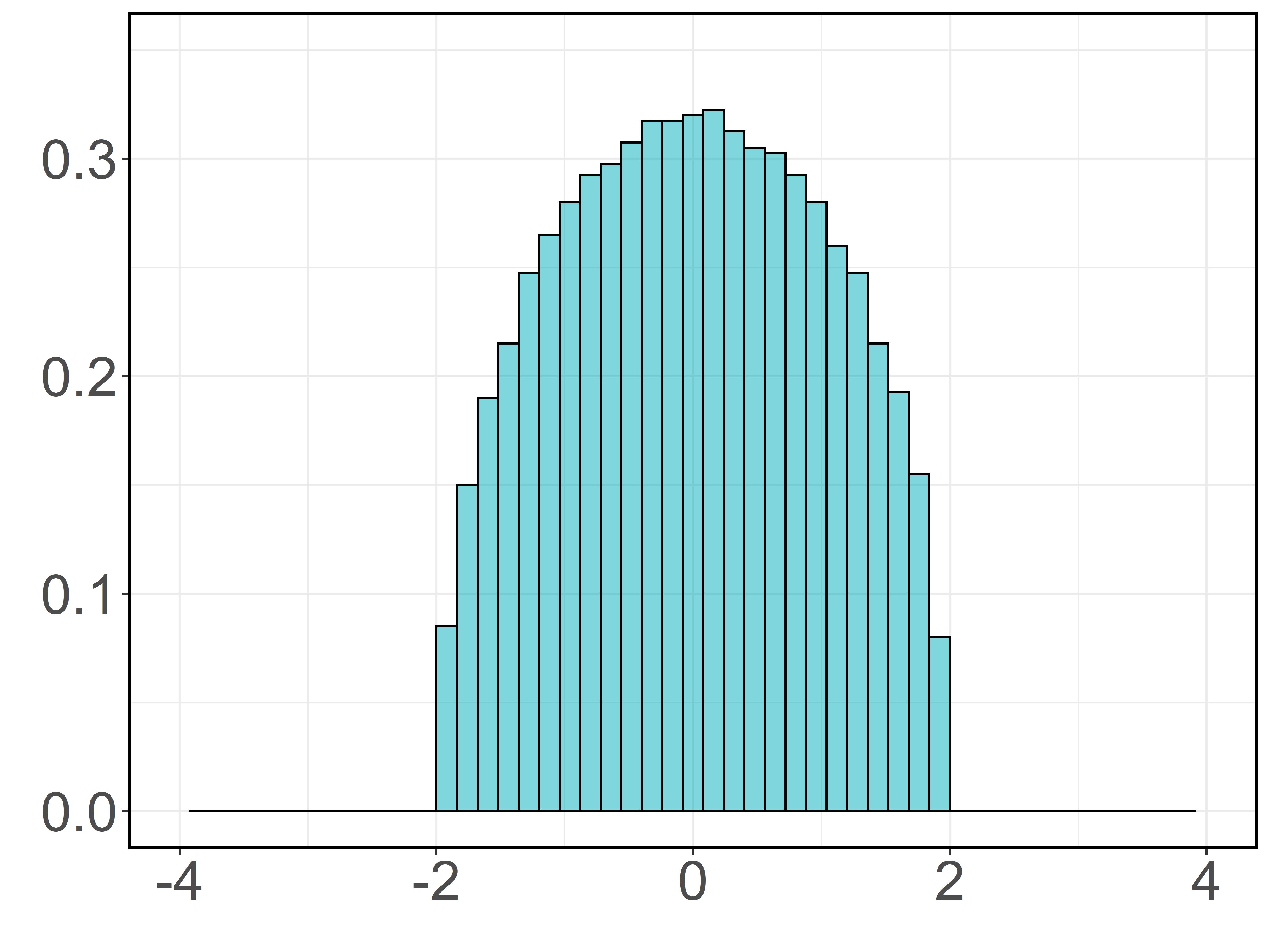}}}%
\subfloat{{\includegraphics[width=7cm]{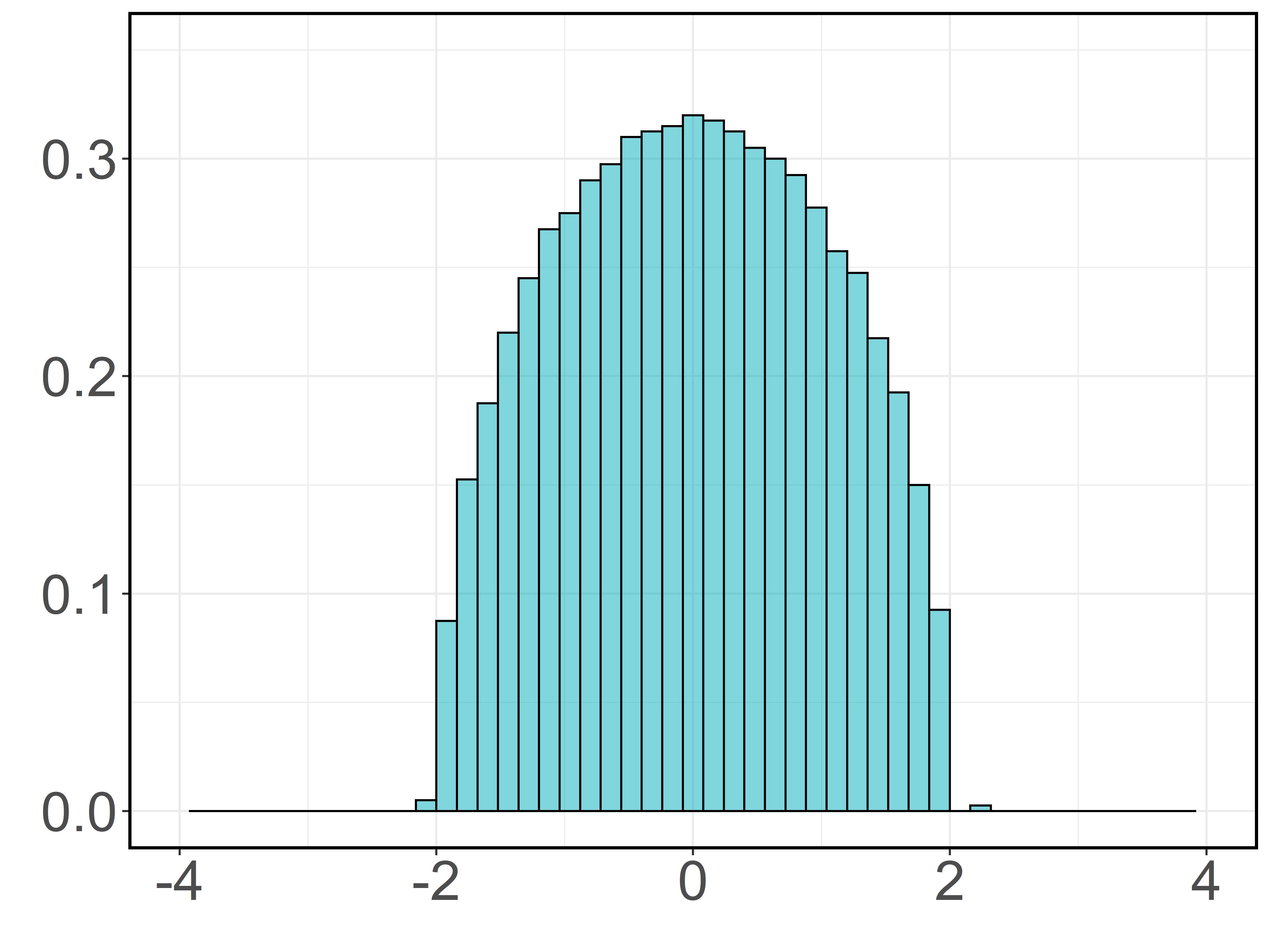}}}%
\caption{The histograms of the eigenvalues for the hidden community model with $p=0.2$ (left) and $p=0.25$ (right).}
\label{fig:hidden}
\end{figure}

We performed the numerical simulation for the hidden community model. We set $N=2500$, $\gamma = 1/4$, and $q=0.2$. Following the analysis in Appendix \ref{subsec:hidden}, we find that an outlier eigenvalue occurs if
\[
	p > q + \frac{\sqrt{q(1-q)}}{\gamma \sqrt{N}} = 0.232.
\]
In Figure \ref{fig:hidden}, we compare the histograms of the eigenvalues of the shifted, rescaled adjacency matrices with $p=0.2$ and $p=0.25$, respectively. As predicted by the analysis, the outlier appears only for the case $p=0.25$.

\subsection{Unbalanced stochastic block model}\label{sec:exu}

We next consider the case $p_1=p_2$ or $\alpha_1=\alpha_2$ with $\gamma \neq 1/2$, which we will refer to an unbalanced stochastic block model. As in the hidden community model, the transition occurs in the regime $p_1=p_2 := p = \frac{w}{\sqrt{N}}+q$. After shifting and rescaling, we find that $\lambda_2 \to 2$ and
\[
	\lambda_1 \to
		\begin{cases}
			\frac{w}{2\sqrt{q(1-q)}} + \frac{2\sqrt{q(1-q)}}{w} & \text{ if } w > 2\sqrt{q(1-q)} \,, \\
			2 & \text{ if } w < 2\sqrt{q(1-q)} \,.
		\end{cases}
\]
See Appendix \ref{subsec:unbal} for the detail. Note that the transition does not depend on $\gamma$.

\begin{figure}[!ht]
\centering
\subfloat{{\includegraphics[width=7cm]{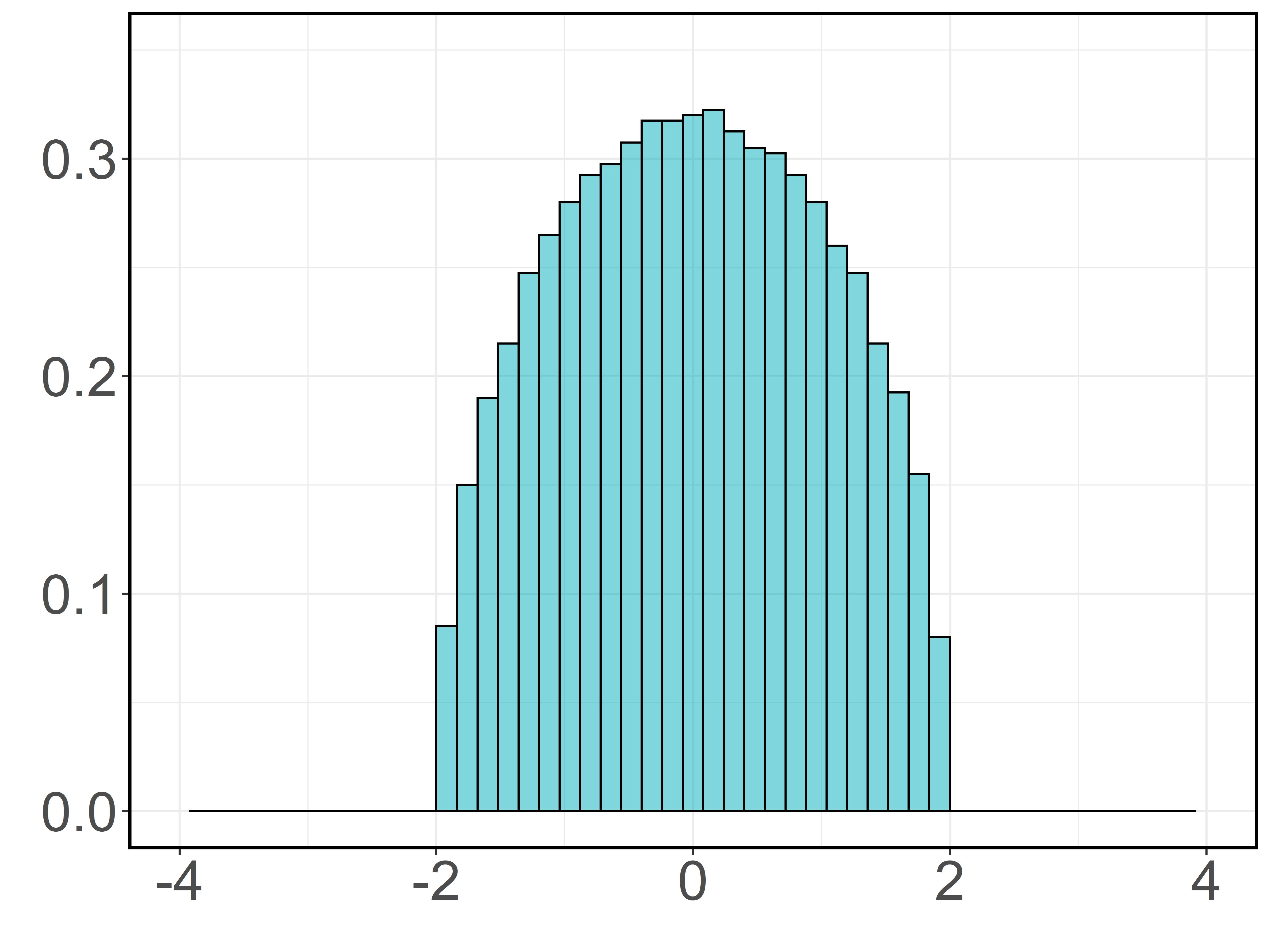}}}%
\subfloat{{\includegraphics[width=7cm]{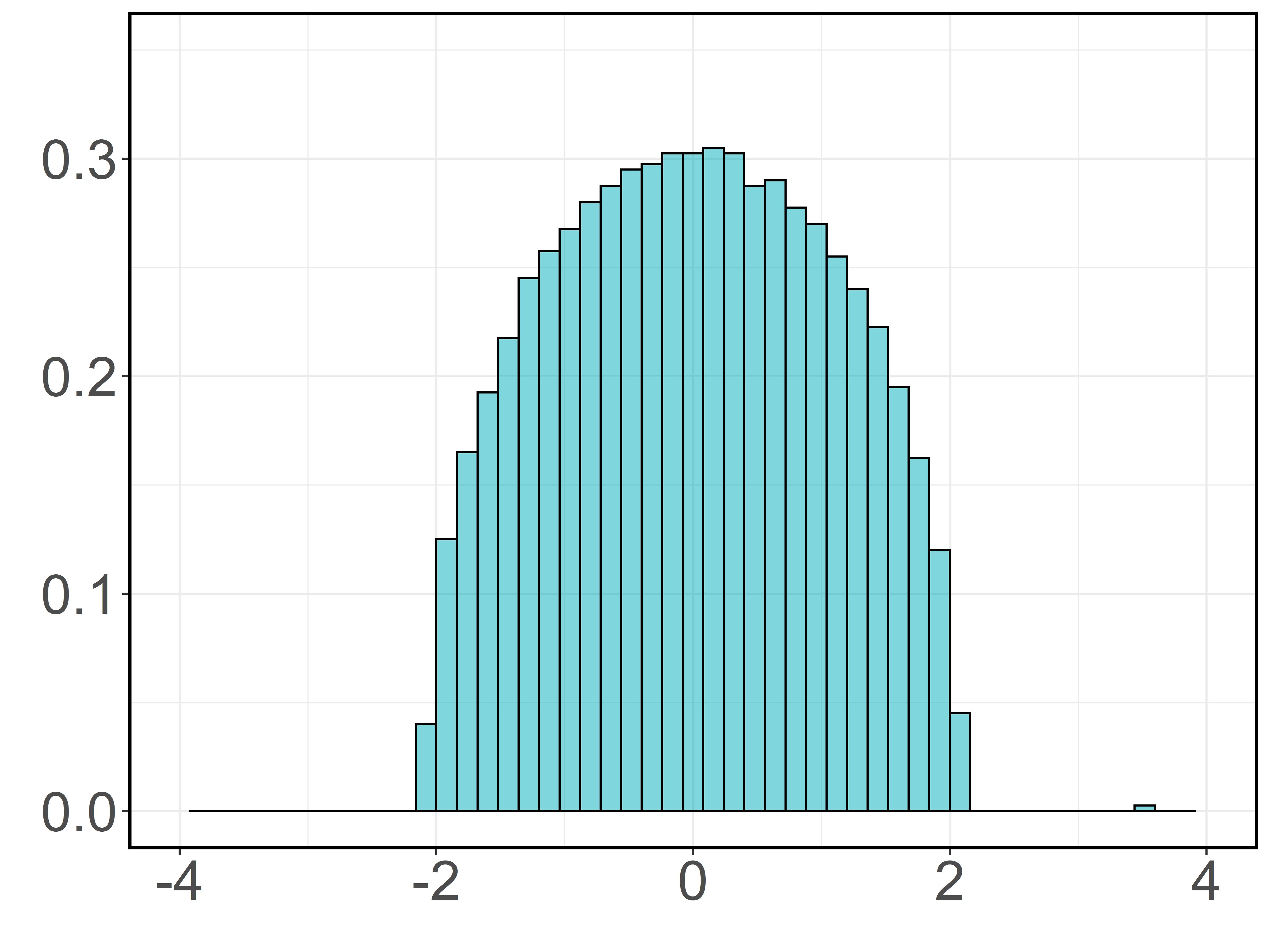}}}%
\caption{The histograms of the eigenvalues for the unbalanced stochastic block model with $p=0.2$ (left) and $p=0.25$ (right).}
\label{fig:unbal}
\end{figure}

We performed the numerical simulation for the unbalanced stochastic block model. As in the hidden community model, we set $N=2500$, $\gamma = 1/4$, and $q=0.2$. An outlier eigenvalue occurs if
\[
	p > q + \frac{2\sqrt{q(1-q)}}{\sqrt{N}} = 0.216.
\]
In Figure \ref{fig:unbal}, we compare the histograms of the eigenvalues of the shifted, rescaled adjacency matrices with $p=0.2$ and $p=0.25$, respectively. Again, as predicted by the analysis, the outlier appears only for the case $p=0.25$.

\section{Proof of Theorem \ref{thm:main}} \label{sec:sketch}

Recall that we denote by $\lambda_1$ and $\lambda_2$ the two largest eigenvalues of $M$. Let $\mu_1$ and $\mu_2$ be the two largest eigenvalues of $H$. From the result on the Wigner-type matrices, we find that $\mu_1$ and $\mu_2$ converge to the rightmost edge of the limiting ESD of $H$.  By the Cauchy interlacing formula, we have the inequality
\[
	\mu_2 \leq \lambda_2 \leq \mu_1 \leq \lambda_1,
\]
which shows that $\lambda_2$ also converges to the rightmost edge of the limiting ESD of $H$. 

From the minimax principle,
\[
	\lambda_1 = \max_{\|x\|=1} \langle x, mx \rangle = \max_{\|x\|=1} \left( \langle x, Hx \rangle + \lambda |\langle x, u \rangle|^2 \right),
\]
which shows that $\lambda_1$ is an increasing function of $\lambda$. Further, since $\lambda_1 \geq \mu_1$ and
\[
	\lambda - \mu_1 \leq \lambda_1 \leq \lambda + \mu_1,
\]
we find that $\lambda_1 - \mu_1 = o(1)$ if $\lambda = o(1)$ and $\lambda_1 > \mu_1 + 1$ if $\lambda > 2\mu_1 + 1$. Thus, since $\lambda_1$ is a continuous function of $\lambda$, conditional on $H$ there exists $\lambda_c$ such that the statement of Theorem \ref{thm:main} holds. It thus remains to show that $\lambda_c$ and $\mu$ are deterministic in the sense that the same statement holds without conditioning on $H$. 

Our proof is based on the Stieltjes transform method in random matrix theory for which we use the following definition:

\begin{defn}[Stieltjes Transform] Let $\mu$ be a probability measure on the real line. The Stieltjes transform of $\mu$ is defined by 
\[
	S_\mu(z)=\int_\mathbb{R}\frac{1}{x-z} \dd\mu(x)
\]
for $z\in \mathbb{C} \backslash \mathrm{supp}\,(\mu)$.
\end{defn}

For the noise $H$, we consider its resolvent $G(z)$ defined by
\beq \label{eq:resolvent}
	G(z) := (H-zI)^{-1}
\eeq
for $z \in \mathbb{C} \backslash \mathrm{spec}\,(H)$. Note that the normalized trace $m:= N^{-1} \Tr G$ is equal to the Stieltjes transform of the empirical spectral distribution (ESD) of $H$.

To find the largest eigenvalue of $M = H + \lambda uu^T$, we recall that any eigenvalue $z$ of $M$ satisfies 
\beq
	\det(H+\lambda uu^T-zI)=0,
\eeq	
which can be further decomposed into
\beq \begin{split} \label{eq:det_eq}
    0 &= \det(H+\lambda uu^T-zI) = \det(H-zI)(I+(H-zI)^{-1}\lambda uu^T) \\
		&= \det(H-zI)\cdot \det(I+(H-zI)^{-1}\lambda uu^T).
\end{split} \eeq
Thus, if $z$ is not an eigenvalue of $H$, we find that $\det(H-zI)=0$, and hence 
\[
	\det(I+(H-zI)^{-1}\lambda uu^T)=0.
\]

We now claim that $(H-zI)^{-1}\lambda uu^T$ has rank one. To prove the claim, we notice that for any $v \in \R^N$,
\[
	(H-zI)^{-1}\lambda uu^T v = \langle u, v \rangle (H-zI)^{-1}\lambda u.
\]
It implies that the range of $(H-zI)^{-1}\lambda uu^T$ is contained in $\mathrm{span}\, ((H-zI)^{-1}\lambda u)$, which is a 1-dimensional space.

Since $(H-zI)^{-1}\lambda uu^T$ has rank 1, it has only one non-zero eigenvalue, which we call $\lambda_0$. Then, $\lambda_0$ is $-1$ for otherwise every eigenvalues of $I+(H-zI)^{-1}\lambda uu^T$ is non-zero, contradicting \eqref{eq:det_eq}. Furthermore, it is also obvious that $(H-zI)^{-1}\lambda u$ is an eigenvector associated with the eigenvalue $-1$. Thus, 
\[
	(H-zI)^{-1}\lambda uu^T(H-zI)^{-1}u=-(H-zI)^{-1}u,
\]
which leads us to the equation
\beq \label{eq:largest_eq}
	u^T(H-zI)^{-1}u = \langle u, G(z) u \rangle = -\frac{1}{\lambda}.
\eeq

For the noise matrix $H$, which is a Wigner-type matrix considered in \cite{AEK17}, we have that
\beq \label{eq:approx_local}
	\langle u, G(z) u \rangle \simeq \sum_{i=1}^N m_i(z) u_i^2
\eeq
for any $z$ not contained in an open neighborhood of the support of the limiting ESD of $H$, where we let $\mathbf m:= (m_1, m_2, \dots, m_N)$ be the solution to the quadratic vector equation (QVE)
\beq \label{eq:QVE}
    -\frac{1}{m_i(z)}=z+ \sum^N_{j=1} \E[H_{ij}^2] m_j(z).
\eeq
is satisfied for $i,\ j=1,2,...,N$. (See Appendix \ref{sec:prelim} for the precise statement of \eqref{eq:approx_local}.) We remark that the uniqueness of the solution $m$ for \eqref{eq:QVE} is also known \cite{ajanki2019quadratic}. 

To solve the equations \eqref{eq:largest_eq} and \eqref{eq:QVE}, we need to estimate $m(z)$ from the assumption on the community structure in Definition \ref{def:GSBM}. Assume for the simplicity that $S = \{ 1, 2, \dots, N_1 \}$. From the symmetry, we have an ansatz
\[
	m_1(z) = m_2(z) = \dots = m_{N_1}(z), \quad m_{N_1 +1}(z) = \dots = m_N(z).
\]
Then, we can rewrite \eqref{eq:QVE} as 
\begin{align*}
-\frac{1}{m_i(z)}=
 \begin{cases}
    z+\sum^{N_1}_{j=1}\frac{\alpha_1}{N}m_j(z)+\sum^{N}_{j=N_1+1}\frac{1}{N}m_j(z) & \text{if } \quad 1\le i\le N_1 \\
    z+\sum^{N_1}_{j=1}\frac{1}{N}m_j(z)+\sum^{N}_{j=N_1+1}\frac{\alpha_2}{N}m_j(z) & \text{if } \quad N_1+1\le i\le N \\
 \end{cases},
\end{align*}
which can be further simplified (after omitting the $z$-dependence) to
\beq \begin{split} \label{eq:system}
     -1&=z m_1+\alpha_1\gamma (m_1) ^2+(1-\gamma) m_1 m_N\, \\
     -1&=z m_N+\gamma m_1m_N+\alpha_2(1-\gamma) (m_N)^2\,.
\end{split} \eeq
We can thus conclude that if there exists real $z$ that solves \eqref{eq:system} under the assumption
\beq \label{eq:assump_crit}
	N_1 m_1 \theta_1^2 + (N-N_1) m_N \theta_2^2 = N \left( \gamma m_1 \theta_1^2 + (1-\gamma) m_N \theta_2^2 \right)  =-\frac{1}{\lambda}
\eeq
then $\lambda_1$ converges to $z$ with high probability. Since \eqref{eq:assump_crit} is deterministic, we find that the gap $g$ is deterministic in the supercritical case $\lambda > \lambda_c$.

It remains to find the critical $\lambda_c$. Recall that we set $m:= N^{-1} \Tr G$. From the Stieltjes inversion formula, the (normalized) imaginary part of $m(z)$ corresponds to the density of the limiting ESD of $H$ at $\re z$. Thus, after changing \eqref{eq:system} as a single equation involving $z$ and $m$ only, i.e., $f(z, m)=0$, we find that the upper edge $L_+$ of the ESD of $H$ is the largest real number such that $f(L_+, m) = 0$ has a double root when considered as an equation for $m$. (Note that technically the condition can be checked by solving $f(L_+, m) = 0$ and $\frac{\partial}{\partial m} f(L_+, m) = 0$ simultaneously.) We can thus conclude that $\lambda_c$ is determined as the largest number such that when $\lambda = \lambda_c$ the solution $z$ for the equation \eqref{eq:system} under the assumption \eqref{eq:assump_crit} coincides with $L_+$. This in particular shows that $\lambda_c$ is also deterministic and completes the proof of Theorem \ref{thm:main}.

\section{Conclusion and Future Works} \label{sec:conclusion}
In this paper, we considered the generalized stochastic block model with two communities. We showed the phase transition in the GSBM where the random part is the Wigner-type matrix, which extends the BBP transition. For the precise formulas, we discussed a hidden community model and unbalanced stochastic block model with Bernoulli distribution and Gaussian distribution at the Kesten--Stigum threshold. Both models can be improved for a non-Gaussian case.

We believe that it is possible to prove the phase transition for the sparse matrix in which the data matrix is not necessarily symmetric and most of the elements are composed of zeros. We also hope to extend our result to the GSBM with more than two communities.

\subsection*{Acknowledgments}
The authors were supported in part by the National Research Foundation of Korea (NRF) grant funded by the Korean government (MSIT) (No. 2019R1A5A1028324).

\begin{appendices}

\section{Local law for Wigner-type matrices} \label{sec:prelim}

In this section, we provide a precise statement of the local law for Wigner-type matrices, which was used in the proof of Theorem \ref{thm:main} in Section \ref{sec:sketch}. Wigner-type matrices are defined as follows:

\begin{defn}(Wigner-type matrix)
We say an $N\times N$ random matrix $H=(H_{ij})$ is a Wigner-type matrix if the entries of $H$ are independent real symmetric variables satisfying the following conditions:
\begin{itemize}
    \item $\E(H_{ij})=0$ for all $i, j$.
    \item The variance matrix $\boldsymbol{S}=(S_{ij})$ where $S_{ij}=\E|H_{ij}|^2$ satisfies \[(S^L)_{ij}\ge \frac{\rho}{N} \text{  and  } s_{ij}\le \frac{S_*}{N},\ \ 1\le i, j\le N\] for finite parameters $\rho, S_*, L$.
\end{itemize}
\end{defn}

For the precise statement of the local law, we use the following definitions, which are frequently used in the analysis involving rare events in random matrix theory.

\begin{defn}(Overwhelming probability)
An event $\Omega$ holds with overwhelming probability if for any big enough $D>0$ , $P(\Omega)\le N^{-D}$ for any sufficient large $N$.
\end{defn}

\begin{defn}(Stochastic domination) Let consider two families of non-negative random variables:
    \[\psi=\{\psi^{(N)}(u)|N\in\N,\ u\in U^{(N)}\}\]
    \[\phi=\{\phi^{(N)}(u)|N\in\N,\ u\in U^{(N)}\}\]
where $U^{(N)}$ is N-dependent parameter set. Suppose $N_0: (0, \infty)^2\to\mathbb{N}$ is a given function depending on $p,\ q,\ n$ and $\mu$. If for $\epsilon>0$ small enough and $D>0$ big enough, we have
\[P\left(\phi^{(N)}>N^\epsilon\psi^{(N)}\right)\le N^{-D}, \ \ \ \text{for }N\ge N_0(\epsilon, D),\] then $\phi$ is stochastically dominated by $\psi$ which denoted by $\phi\prec\psi$.
\end{defn}

We are now ready to state the local law. Let $m_i(z)$ be the solution of QVE in \eqref{eq:QVE} and $\rho$ is the density defined as
\[
	\rho(\tau) := \lim_{\rho \searrow 0} \frac{1}{\pi N} \sum_{i=1}^N \im m_i (\tau+ \ii \eta).
\]
(See also Corollary 1.3 of \cite{AEK17} for more detail.)

\begin{thm}[Local law] \label{thm:local_law} Let $H$ be a Wigner-type matrix and fix an arbitrary $\gamma\in(0,1)$. Then, uniformly for all $z=a+bi$ with $b\ge N^{\gamma-1}$, the resolvent $G(z) = (H-zI)^{-1}$ satisfy
\[max_{i,j}\left\vert G_{ij}(z)- m_i(z)\delta_{ij}\right\vert\prec \frac{1+\sqrt{\rho(z)}}{\sqrt{bN}}+\frac{1}{bN}\]
Furthermore, for any deterministic vector $w\in \C^N$ with $\max_i |w_i|\ge 1$, we have
\[\left\vert\sum^N_{i, j=1} \overline{w_i}\left(G_{ij}(z)- m_i(z)\right)\right\vert\prec \frac{1}{\sqrt{bN}}\]
\end{thm}

The local law can be generalized to the anisotropic local law as follows.

\begin{thm}[Anisotropic law] Suppose that the assumptions in Theorem \ref{thm:local_law} hold. Then, uniformly for all $z=a+bi$ with $b\ge N^{\gamma-1}$, and for any two deterministic $\ell^2$-normalized vectors $w, v \in \C^N$,  we have
\[\left\vert\sum^N_{i, j=1} \overline{w_i}G_{ij}(z)v_j-\sum^N_{i=1} m_i(z)\overline{w_i}v_j\right\vert\prec \frac{1+\sqrt{\rho(z)}}{\sqrt{bN}}+\frac{1}{bN}\]
\end{thm}

\section{Examples from stochastic block models} \label{sec:sbm}

In this appendix, we consider stochastic block models, which corresponds to GSBMs with Bernoulli distribution in our setting. Suppose that $\widehat{H}=[\widehat{H}_{ij}]_{i, j=1}^N$ is an SBM such that
\beq \label{eq:H_def}
\widehat{H}=
 \begin{cases}
    \widehat{H}_{ij}\sim Bernoulli(p_1), & \mbox{ if } 1\le i, j\le N_1\,, \\
    \widehat{H}_{ij}\sim Bernoulli(p_2), & \mbox{ if } N_1+1\le i, j\le N\,, \\
    \widehat{H}_{ij}\sim Bernoulli(q), & \mbox{ otherwise}\,. \\
 \end{cases}
\eeq
In what follows, we will call the $(i, j)$-entry is in the diagonal block if $1\le i, j\le N_1$ or $N_1+1\le i, j\le N$, and  otherwise it is in the off-diagonal block.
In the block matrix form, it can also be expressed as follows:
\begin{align*}
\newcommand\explainA{%
  \overbrace{%
    \hphantom{\begin{matrix}Bernoulli(p)\end{matrix}}%
  }^{\text{$N_1$}}%
}
\newcommand{\explainB}{%
  \overbrace{%
    \hphantom{\begin{matrix}Bernoulli(p)\end{matrix}}%
  }^{\text{$N-N_1$}}%
}
\newcommand{\explainC}{%
  \left.\vphantom{\begin{matrix}Bernoulli(p)\end{matrix}}\right\}%
  \text{\scriptsize$N_1$}%
}
\newcommand{\explainD}{%
  \left.\vphantom{\begin{matrix}Bernoulli(p)\end{matrix}}\right\}%
  \text{\scriptsize$N-N_1$}%
}
\widehat{H}=
\settowidth{\dimen0}{%
  $\begin{pmatrix}\vphantom{\begin{matrix}\lambda_1&\dots&\lambda_1\end{matrix}}\end{pmatrix}$%
}
\settowidth{\dimen2}{$\explainB$}
\begin{matrix}
\begin{matrix}\hspace*{0.5\dimen0}\explainA&\explainB\hspace*{0.5\dimen0}\end{matrix}
\\[-0.5ex]
\begin{pmatrix}
\begin{matrix}
Bernoulli(p_1)
\end{matrix}
&
\begin{matrix}
Bernoulli(q)
\end{matrix}
&\\
\begin{matrix}
Bernoulli(q)
\end{matrix}
&
\begin{matrix}
Bernoulli(p_2)
\end{matrix}
\end{pmatrix}
&\hspace*{-1em}\begin{matrix} \explainC\hfill \\ \explainD \end{matrix}
\end{matrix}
\end{align*}

Our goal is to shift and rescale $\wh H$ to convert it into a GSBM $M = H + \lambda uu^T$ in Definition \ref{def:GSBM}. We first notice that the variances of the entries of $\wh H$ are $p_1 (1-p_1)$ and $p_2 (1-p_2)$ for the diagonal block and $q(1-q)$ for the off-diagonal block. Since we assume that the variance of the entry $H_{ij}$ in the off-diagonal block is $N^{-1}$, we find that the matrix must be divided by $\sqrt{Nq(1-q)}$. It is then immediate to find that
\beq
	\alpha_1=\frac{p_1(1-p_1)}{q(1-q)}, \qquad \alpha_2=\frac{p_2(1-p_2)}{q(1-q)}.
\eeq
as in \eqref{eq:p_12}.

The mean matrix
\[
	\E[\wh H] =
		\begin{pmatrix}
		p_1 & q \\
		q & p_2
		\end{pmatrix}
\]
is a rank-$2$ matrix, and thus we need to subtract each entry by a deterministic number, which depends on the parameters $p_1, p_2$, and $q$.

\subsection{Hidden community model} \label{subsec:hidden}
Suppose that $p_1 = p$ and $p_2 = q$. It is then easy to find that $\E[\wh H]$ becomes a rank-$1$ matrix after subtracting each entry by $q$, i.e., if we let $E_0$ be the $N \times N$ matrix whose all entries are $q$, then
\[
	\E[\wh H] - E_0 =
		\begin{pmatrix}
		p-q & 0 \\
		0 & 0
		\end{pmatrix}.
\]
Thus, we find that
\beq
	M = \frac{1}{\sqrt{Nq(1-q)}} (\wh H - E_0)
\eeq
and
\beq
	\E[M] = \frac{1}{\sqrt{Nq(1-q)}}
		\begin{pmatrix}
		p-q & 0 \\
		0 & 0
		\end{pmatrix}.
\eeq
Recall that $N_1 = \gamma N$ and $p = \frac{w}{\sqrt{N}} + q$. Since $\lambda u u^T = \E[M]$, we get
\begin{align*}
\newcommand\explainA{%
  \overbrace{%
    \hphantom{\begin{matrix}-\gamma(1-\gamma)(p_1+p_2-2q)\end{matrix}}%
  }^{\text{$N_1$}}%
}
\newcommand{\explainB}{%
  \overbrace{%
    \hphantom{\begin{matrix}-\gamma(1-\gamma)(p_1+p_2-2q)\end{matrix}}%
  }^{\text{$N-N_1$}}}
\newcommand{\explainC}{%
  \left.\vphantom{\begin{matrix}-\gamma(1-\gamma)(p_1+p_2-2q)\end{matrix}}\right\}%
  \text{\scriptsize$N_1$}%
}
\newcommand{\explainD}{%
  \left.\vphantom{\begin{matrix}-\gamma(1-\gamma)(p_1+p_2-2q)\end{matrix}}\right\}%
  \text{\scriptsize$N-N_1$}%
}
\settowidth{\dimen0}{%
  $\begin{pmatrix}\vphantom{\begin{matrix}\lambda_1&\dots&\lambda_1\end{matrix}}\end{pmatrix}$%
}
u=
\settowidth{\dimen2}{$\explainB$}
\begin{matrix}
\\[-1.5em]
\begin{pmatrix}
\begin{matrix}
\frac{1}{\sqrt{\gamma N}}
\end{matrix}
&\\[0.5em]
\begin{matrix}
0
\end{matrix}
\end{pmatrix}
&\hspace*{-1em}\begin{matrix} \explainC\hfill \\ \explainD \end{matrix}
\end{matrix},
\end{align*}
i.e., $\theta_1 = 1/\sqrt{\gamma N}$ and $\theta_2 = 0$. We also find that
\[
	\lambda = \frac{N_1 (p-q)}{\sqrt{Nq(1-q)}} = \frac{\gamma w}{\sqrt{q(1-q)}}.
\]

Following the proof of Theorem \ref{thm:main} in Section \ref{sec:sketch}, we solve the system of equation in \eqref{eq:system},
\beq \begin{split} \label{eq:system1}
     -1&=z m_1+ \frac{p(1-p)}{q(1-q)} \gamma (m_1)^2+(1-\gamma) m_1 m_N \,, \\
     -1&=z m_N+\gamma m_1 m_N+ (1-\gamma) (m_N)^2 \,.
\end{split} \eeq
Since $p-q = O(N^{-1/2})$, we consider an ansatz $m_N = m_1 + O(N^{-1/2})$, which shows for $m = \gamma m_1 + (1-\gamma)m_N$ that
\[
	1 + zm + m^2 = O(N^{-1/2}).
\]
Following the analysis in the last paragraph of Section \ref{sec:sketch}, we find that the upper edge $L_+ = 2 + O(N^{-1/2})$. By Theorem \ref{thm:main}, it also implies that $\lambda_2 \to 2$ as $N \to \infty$.

In order to determine the location of the largest eigenvalue $\lambda_1$, we consider \eqref{eq:system1} under the assumption in \eqref{eq:assump_crit},
\beq \label{eq:m1h}
	N \left( \gamma m_1 \theta_1^2 + (1-\gamma) m_N \theta_2^2 \right) = m_1 = -\frac{1}{\lambda} = -\frac{\sqrt{q(1-q)}}{\gamma w}.
\eeq
We remark that the ansatz $m_N = m_1 + O(N^{-1/2})$ can be directly checked in this case; by plugging \eqref{eq:m1h} into \eqref{eq:system1} and eliminating $z$,
\[
	\left( \frac{p(1-p)}{q(1-q)} - 1 \right) \gamma \left( \frac{\sqrt{q(1-q)}}{\gamma w} \right)^2 m_N + m_N + \frac{\sqrt{q(1-q)}}{\gamma w} = 0,
\]
whose solution is
\[
	m_N = -\frac{\sqrt{q(1-q)}/(\gamma w)}{1 + \left( \frac{p(1-p)}{q(1-q)} - 1 \right) \gamma \left( \frac{\sqrt{q(1-q)}}{\gamma w} \right)^2} = -\frac{\sqrt{q(1-q)}}{\gamma w} \left( 1-\frac{\sqrt{N} (1-2 q)-w}{N \gamma w+\sqrt{N} (1-2 q)-w} \right).
\]
To find the location of the largest eigenvalue, we need to check whether the assumption \eqref{eq:m1h} is valid. However, we can instead find the value of $z$ by first assuming that the solution exists. Then,
\[
	z = \frac{\gamma w}{\sqrt{q(1-q)}} + \frac{\sqrt{q(1-q)}}{\gamma w} + O(N^{-1/2}).
\]
At the critical $\lambda_c$ for the phase transition in Theorem \ref{thm:main}, the location of the largest eigenvalue coincides with the location of the upper edge $L_+$ in the limit $N \to \infty$, or equivalently, $\frac{\gamma w}{\sqrt{q(1-q)}} = 1$. Thus, we conclude that
\[
	\lambda_1 \to
		\begin{cases}
			\frac{\gamma w}{\sqrt{q(1-q)}} + \frac{\sqrt{q(1-q)}}{\gamma w} & \text{ if } w > \frac{\sqrt{q(1-q)}}{\gamma} \,, \\
			2 & \text{ if } w < \frac{\sqrt{q(1-q)}}{\gamma} \,.
		\end{cases}
\]

\subsection{Unbalanced stochastic model} \label{subsec:unbal}
Suppose that $p_1 = p_2 = p$. Following the strategy in Appendix \ref{subsec:hidden}, we let $E_1$ be the $N \times N$ matrix whose all entries are $(p+q)/2$. Then,
\[
	\E[\wh H] - E_1 =
		\begin{pmatrix}
		(p-q)/2 & (q-p)/2 \\
		(q-p)/2 & (p-q)/2
		\end{pmatrix}.
\]
Thus, we find that
\beq
	M = \frac{1}{\sqrt{Nq(1-q)}} (\wh H - E_1)
\eeq
and
\beq
	\E[M] = \frac{1}{\sqrt{Nq(1-q)}}
		\begin{pmatrix}
		(p-q)/2 & (q-p)/2 \\
		(q-p)/2 & (p-q)/2
		\end{pmatrix}.
\eeq
From $\lambda u u^T = \E[M]$, we get
\begin{align*}
\newcommand\explainA{%
  \overbrace{%
    \hphantom{\begin{matrix}-\gamma(1-\gamma)(p_1+p_2-2q)\end{matrix}}%
  }^{\text{$N_1$}}%
}
\newcommand{\explainB}{%
  \overbrace{%
    \hphantom{\begin{matrix}-\gamma(1-\gamma)(p_1+p_2-2q)\end{matrix}}%
  }^{\text{$N-N_1$}}}
\newcommand{\explainC}{%
  \left.\vphantom{\begin{matrix}-\gamma(1-\gamma)(p_1+p_2-2q)\end{matrix}}\right\}%
  \text{\scriptsize$N_1$}%
}
\newcommand{\explainD}{%
  \left.\vphantom{\begin{matrix}-\gamma(1-\gamma)(p_1+p_2-2q)\end{matrix}}\right\}%
  \text{\scriptsize$N-N_1$}%
}
\settowidth{\dimen0}{%
  $\begin{pmatrix}\vphantom{\begin{matrix}\lambda_1&\dots&\lambda_1\end{matrix}}\end{pmatrix}$%
}
u=
\settowidth{\dimen2}{$\explainB$}
\begin{matrix}
\\[-1.5em]
\begin{pmatrix}
\begin{matrix}
\frac{1}{\sqrt{N}}
\end{matrix}
&\\[0.5em]
\begin{matrix}
-\frac{1}{\sqrt{N}}
\end{matrix}
\end{pmatrix}
&\hspace*{-1em}\begin{matrix} \explainC\hfill \\ \explainD \end{matrix}
\end{matrix},
\end{align*}
i.e., $\theta_1 = 1/\sqrt{N}$ and $\theta_2 = -1/\sqrt{N}$. Also,
\[
	\lambda = \frac{N(p-q)}{2\sqrt{Nq(1-q)}} = \frac{w}{2\sqrt{q(1-q)}}.
\]

With $\alpha_1 = \alpha_2 = \frac{p(1-p)}{q(1-q)}$, we solve the system of equation in \eqref{eq:system},
\beq \begin{split} \label{eq:system2}
     -1&=z m_1+ \frac{p(1-p)}{q(1-q)} \gamma (m_1)^2+(1-\gamma) m_1 m_N \,, \\
     -1&=z m_N+\gamma m_1 m_N+ \frac{p(1-p)}{q(1-q)} (1-\gamma) (m_N)^2 \,.
\end{split} \eeq
Again, we consider an ansatz $m_N = m_1 + O(N^{-1/2})$, which leads us to the result that the upper edge $L_+ = 2 + O(N^{-1/2})$ and $\lambda_2 \to 2$ as $N \to \infty$. The assumption in \eqref{eq:assump_crit} becomes
\beq \label{eq:m_u}
	\gamma m_1 + (1-\gamma) m_N = -\frac{1}{\lambda} = -\frac{2\sqrt{q(1-q)}}{w}.
\eeq
If the solution to the equation \eqref{eq:system2} exists, it would be
\[
	z = \frac{w}{2\sqrt{q(1-q)}} + \frac{2\sqrt{q(1-q)}}{w} + O(N^{-1/2}).
\]
At the critical $\lambda_c$, $\frac{w}{2\sqrt{q(1-q)}} = 1$, and thus we conclude that
\[
	\lambda_1 \to
		\begin{cases}
			\frac{w}{2\sqrt{q(1-q)}} + \frac{2\sqrt{q(1-q)}}{w} & \text{ if } w > 2\sqrt{q(1-q)} \,, \\
			2 & \text{ if } w < 2\sqrt{q(1-q)} \,.
		\end{cases}
\]

\end{appendices}


\begin{thebibliography}{10}

\bibitem{AbbeEmm2017}
E.~Abbe.
\newblock Community detection and stochastic block models: recent developments.
\newblock {\em The Journal of Machine Learning Research}, 18(1):6446--6531, 2017.

\bibitem{AbbeEmm2014}
E.~Abbe, A.~Bandeira, and G.~Hall.
\newblock Exact recovery in the stochastic block model.
\newblock {\em IEEE Transactions on Information Theory}, 62(1):471--487, 2014.

\bibitem{ajanki2019quadratic}
O.~Ajanki, L.~Erd{\H{o}}s, and T.~Kr{\"u}ger.
\newblock Quadratic vector equations on complex upper half-plane.
\newblock{\em American Mathematical Society}, 261(1261), 2019.

\bibitem{AEK17}
O.~Ajanki, L.~Erd{\H{o}}s, and T.~Kr{\"u}ger.
\newblock Universality for general wigner-type matrices.
\newblock {\em Probability Theory and Related Fields}, 169(3):667--727, 2017.

\bibitem{BBP2005}
J.~Baik, G.~Ben~Arous, and S.~P\'ech\'e.
\newblock Phase transition of the largest eigenvalue for nonnull complex sample
  covariance matrices.
\newblock {\em The Annals of Probability}, 33(5):1643--1697, 2005.

\bibitem{NIPS_Barbier}
J.~Barbier, M.~Dia, N.~Macris, F.~Krzakala, T.~Lesieur, and L.~Zdeborov\'{a}.
\newblock Mutual information for symmetric rank-one matrix estimation: A proof
  of the replica formula.
\newblock {\em Advances in Neural Information Processing Systems}, 29:424–432, 2016.

\bibitem{Raj2011}
F.~Benaych-Georges and R.~R. Nadakuditi.
\newblock The eigenvalues and eigenvectors of finite, low rank perturbations of
  large random matrices.
\newblock {\em Advances in Mathematics}, 227(1):494--521, 2011.

\bibitem{CapitaineDonatiFeral2009}
M.~Capitaine, C.~Donati-Martin, and D.~F\'eral.
\newblock The largest eigenvalues of finite rank deformation of large {W}igner
  matrices: convergence and nonuniversality of the fluctuations.
\newblock {\em The Annals of Probability}, 37(1):1--47, 2009.

\bibitem{Chen-2015}
P.-Y. Chen and A.~Hero.
\newblock Universal phase transition in community detectability under a
  stochastic block model.
\newblock {\em Physical Review E}, 91(3):032804, 2015.

\bibitem{chung2019}
H.~W. Chung and J.~O. Lee.
\newblock Weak detection of signal in the spiked wigner model.
\newblock In {\em International Conference on Machine Learning}, 97:1233--1241, 2019.

\bibitem{dumitriu2019sparse}
I.~Dumitriu and Y.~Zhu.
\newblock Sparse general {W}igner-type matrices: Local law and eigenvector
  delocalization.
\newblock {\em Journal of Mathematical Physics}, 60(2):023301, 2019.

\bibitem{AlaouiJordan2018}
A.~El~Alaoui, F.~Krzakala, and M.~I. Jordan.
\newblock {Fundamental limits of detection in the spiked {W}igner model}.
\newblock {\em The Annals of Statistics}, 48(2):863--885, 2020.

\bibitem{erdHos2020cusp}
L.~Erd{\H{o}}s, T.~Kr{\"u}ger, and D.~Schr{\"o}der.
\newblock Cusp universality for random matrices {I}: local law and the complex
  hermitian case.
\newblock {\em Communications in Mathematical Physics}, 378(2):1203--1278,
  2020.

\bibitem{erdHos2019bounds}
L.~Erd{\H{o}}s and P.~M{\"u}hlbacher.
\newblock Bounds on the norm of {W}igner-type random matrices.
\newblock {\em Random Matrices: Theory and Applications}, 8(03):1950009, 2019.

\bibitem{FeralPeche2007}
D.~F\'eral and S.~P\'ech\'e.
\newblock The largest eigenvalue of rank one deformation of large {W}igner
  matrices.
\newblock {\em Communications in mathematical physics}, 272(1):185--228, 2007.

\bibitem{Hajek-2017}
B.~Hajek, Y.~Wu, and J.~Xu.
\newblock Information limits for recovering a hidden community.
\newblock {\em IEEE Transactions on Information Theory}, 63(8):4729-4745, 2017.

\bibitem{hajek2016achieving}
B.~Hajek, Y.~Wu, and J.~Xu.
\newblock Achieving exact cluster recovery threshold via semidefinite
  programming.
\newblock {\em IEEE Transactions on Information Theory}, 62(5):2788--2797,
  2016.

\bibitem{hajek2018recovering}
B.~Hajek, Y.~Wu, and J.~Xu.
\newblock Recovering a hidden community beyond the {K}esten--{S}tigum threshold in
  $O(|E| log*|V|)$ time.
\newblock {\em Journal of Applied Probability}, 55(2):325--352, 2018.

\bibitem{jung2020}
J.~H. Jung, H.~W. Chung, and J.~O. Lee.
\newblock Weak detection in the spiked wigner model with general rank.
\newblock {\em arXiv:2001.05676}, 2020.

\bibitem{Peche2006}
S.~P\'ech\'e.
\newblock The largest eigenvalue of small rank perturbations of {H}ermitian
  random matrices.
\newblock {\em Probability Theory and Related Fields}, 134(1):127--173, 2006.

\bibitem{Perry2018}
A.~Perry, A.~S. Wein, A.~S. Bandeira, and A.~Moitra.
\newblock Optimality and sub-optimality of {PCA} {I}: {S}piked random matrix
  models.
\newblock {\em The Annals of Statistics}, 46(5):2416--2451, 2018.

\bibitem{stanley2019stochastic}
N.~Stanley, T.~Bonacci, R.~Kwitt, M.~Niethammer, and P.~J. Mucha.
\newblock Stochastic block models with multiple continuous attributes.
\newblock {\em Applied Network Science}, 4(1):1--22, 2019.

\bibitem{zhu2020graphon}
Y.~Zhu.
\newblock A graphon approach to limiting spectral distributions of wigner-type
  matrices.
\newblock {\em Random Structures \& Algorithms}, 56(1):251--279, 2020.

\end{thebibliography}
\end{document}